\documentclass[9pt]{article}
\usepackage{mathrsfs}
\usepackage{amsthm}
\usepackage{amssymb}
\usepackage{amsmath}
\usepackage{graphicx}
\usepackage{color}
\usepackage{amsfonts}
\usepackage{float}
\usepackage{cite}
\usepackage [latin1]{inputenc}
\usepackage[text={140mm,210mm},left=35mm,vmarginratio=1:1]{geometry}
\newtheorem{theorem}{Theorem}[section]
\newtheorem{remark}{Remark}[section]
\newtheorem{lemma}[theorem]{Lemma}
\newtheorem{proposition}[theorem]{Proposition}

\numberwithin{equation}{section}
\normalsize

\begin{document}
\title{\textbf{Moderate deviations of generalized $N$-urn Ehrenfest models}}

\author{Lirong Ren \thanks{\textbf{E-mail}: 20121634@bjtu.edu.cn \textbf{Address}: School of Science, Beijing Jiaotong University, Beijing 100044, China.}\\ Beijing Jiaotong University\\ Xiaofeng Xue \thanks{\textbf{E-mail}: xfxue@bjtu.edu.cn \textbf{Address}: School of Science, Beijing Jiaotong University, Beijing 100044, China.}\\ Beijing Jiaotong University}

\date{}
\maketitle

\noindent {\bf Abstract:} This paper is a further investigation of the generalized $N$-urn Ehrenfest model introduced in \cite{Xue2020}. A moderate deviation principle from the hydrodynamic limit of the model is derived. The proof of this main result follows a routine procedure introduced in \cite{Kipnis1989}, where a replacement lemma plays the key role. To prove the replacement lemma, the large deviation principle of the model given in \cite{Xue2020} is utilized.

\quad

\noindent {\bf Keywords:} hydrodynamic limit, N-urn Ehrenfest model, moderate deviation, replacement lemma.

\section{Introduction and main results}\label{section one}
In this paper we will prove a moderate deviation principle from the hydrodynamic limit of the generalized $N$-urn Ehrenfest model introduced in \cite{Xue2020}. We first recall the definition of the model. Initially some gas molecules are put into $N$ boxes, where $N\geq 2$ is an integer. We assume that numbers of gas molecules in different boxes are independent and the number of gas molecules in the $i$th box follows Poisson distribution with mean $\phi(\frac{i}{N})$ for $1\leq i\leq N$, where $\phi$ is a positive function in $C\left([0, 1]\right)$. For $1\leq i,j\leq N$, each gas molecule in the $i$th box jumps to the $j$th box at rate $\frac{1}{N}\lambda\left(\frac{i}{N}, \frac{j}{N}\right)$, where $\lambda$ is a positive function in $C^{1,1}\left([0, 1]\times[0, 1]\right)$. When $\lambda\equiv 1$, the above model reduces to the classic $N$-urn Ehrenfest model introduced in \cite{Cheng2020}. For any $t\geq 0$, let $X_t^N(i)$ be the number of gas molecules in the $i$th box at moment $t$ and
\[
X_t^N=\left(X_t^N(1), X_t^N(2),\ldots, X_t^N(N)\right),
\]
then $\{X_t^N\}_{t\geq 0}$ is a continuous-time Markov process with state space $\{0, 1,2,\ldots\}^N$ and generator $\mathcal{L}_N$ given by
\[
\mathcal{L}_Nf(x)=\sum_{i=1}^N\sum_{j=1}^N\frac{x(i)}{N}\lambda\left(\frac{i}{N}, \frac{j}{N}\right)\left[f(x^{i,j})-f(x)\right]
\]
for any $f\in C\left(\{0, 1,2,\ldots\}^N\right)$ and $x\in \{0, 1,2,\ldots\}^N$, where $x^{i,j}=x$ when $i=j$ and
\[
x^{i,j}(l)=
\begin{cases}
x(l) & \text{~if~}l\neq i,j,\\
x(i)-1 & \text{~if~}l=i,\\
x(j)+1 & \text{~if~}l=j
\end{cases}
\]
when $i\neq j$. 

Now we recall the hydrodynamic limit of $\{X_t^N\}_{t\geq 1}$ given in \cite{Xue2020}. For each $N\geq 1$ and any $t\geq 0$, we define the empirical measure $\mu_t^N$ as
\[
\mu_t^N(du)=\frac{1}{N}\sum_{i=1}^NX_t^N(i)\delta_{\frac{i}{N}}(du),
\]
where $\delta_{\frac{i}{N}}(du)$ is the Dirac measure concentrated at $\frac{i}{N}$. That is to say, $\mu_t^N$ is a random linear operator from $C([0, 1])$ to $\mathbb{R}$ that
\[
\mu_t^N(f)=\int_{[0, 1]} f(u)\mu_t^N(du)=\frac{1}{N}\sum_{i=1}^NX_t^N(i)f(\frac{i}{N})
\]
for any $f\in C([0, 1])$. Let $P_1$ be the linear operator from $C([0, 1])$ to $C([0, 1])$ that
\[
(P_1f)(x)=\int_0^1\lambda(x,y)f(y)dy
\]
for any $f\in C([0,1]), x\in [0,1]$ and $P_2$ be the one that
\[
(P_2f)(x)=\int_0^1\lambda(x,y)f(x)dy
\]
for any $f\in C([0, 1]), x\in [0,1]$, then it is shown in \cite{Xue2020} that there is a unique deterministic measure-valued process $\{\mu_t\}_{t\geq 0}$ that
\begin{equation}\label{equation HydrodynamicLimit}
\mu_t(f)=\int_0^1f(x)\phi(x)dx+\int_0^t\mu_s\left((P_1-P_2)f\right)ds
\end{equation}
for any $t\geq 0$ and $f\in C([0, 1])$. The following proposition is proved in \cite{Xue2020}, which gives the hydrodynamic limit of $\{X_t^N\}_{t\geq 0}$ as $N\rightarrow+\infty$.
\begin{proposition}[{\cite[Theorem 2.3]{Xue2020}}]\label{proposition 1.1 HydrodynamicLimit}
 Let $\mu$ be defined as in Equation \eqref{equation HydrodynamicLimit}, then
\[
\lim_{N\rightarrow+\infty}\mu_t^N(f)=\mu_t(f)
\]
in probability for any $t\geq 0$ and $f\in C([0, 1])$.
\end{proposition}

In this paper, we are concerned with the moderate deviation principle from the hydrodynamic limit given in Proposition \ref{proposition 1.1 HydrodynamicLimit}. To give our results, we first introduce some notations and definitions. We use $\mathcal{S}$ to denote the dual of $C([0, 1])$, i.e., the set of linear operators from $C([0, 1])$ to $\mathbb{R}$. For later use, we use $\mathcal{A}$ to denote the subset of $\mathcal{S}$ consist of nonnegative measures, i.e,
$\nu\in \mathcal{A}$ if and only if $\nu(f)\geq 0$ for any nonnegative $f\in C([0, 1])$. For given $T_0>0$, we use $\mathcal{D}([0, T_0], \mathcal{S})$ to denote the set of c\`{a}dl\`{a}g functions from $[0, T_0]$ to $\mathcal{S}$. For any $\nu\in \mathcal{S}$, we define
\begin{equation}\label{equ 1.2 definition of Iini}
I_{ini}(\nu)=\sup_{f\in C([0, 1])}\left\{\nu(f)-\frac{1}{2}\int_0^1\phi(x)f^2(x)dx\right\}.
\end{equation}
For any $\pi\in \mathcal{D}([0, T_0], \mathcal{S})$, we define
\begin{align}\label{equ 1.3 defintion of Idyn}
I_{dyn}(\pi)&=\sup_{G\in C^{1,0}([0, T_0]\times [0, 1])}\Bigg\{\pi_{T_0}(G_{T_0})-\pi_0(G_0)-\int_0^{T_0}\pi_s\left((\partial_s+P_1-P_2)G_s\right)ds \notag\\
&-\frac{1}{2}\int_0^{T_0}\left(\int_{[0,1]}\left(\int_0^1\lambda(x,y)\left(G_s(y)-G_s(x)\right)^2dy\right)\mu_s(dx)\right)ds\Bigg\},
\end{align}
where $\mu$ is defined as in Equation \eqref{equation HydrodynamicLimit} and $G_t(\cdot)=G(t, \cdot)$ for any $G\in C^{1,0}([0, T_0]\times [0, 1]), 0\leq t\leq T_0$. Let $\{a_N\}_{N\geq 1}$ be a given positive sequence that $\lim_{N\rightarrow+\infty}\frac{a_N}{N}=\lim_{N\rightarrow+\infty}\frac{\sqrt{N}}{a_N}=0$, then we define random measure $\theta_t^N$ as
\[
\theta_t^N(du)=\frac{1}{a_N}\sum_{i=1}^N\left(X_t^N(i)-EX_t^N(i)\right)\delta_{\frac{i}{N}}(du)
\]
for any $N\geq 1$ and $0\leq t\leq N$. We use $\theta^N$ to denote $\{\theta_t^N:~0\leq t\leq T_0\}$, then $\theta^N\in \mathcal{D}([0, T_0], \mathcal{S})$. Now we give our main result.

\begin{theorem}\label{theorem 1.1 MDP}Let $I_{ini}$ and $I_{dyn}$ be defined as in Equations \eqref{equ 1.2 definition of Iini} and \eqref{equ 1.3 defintion of Idyn} respectively, then
\begin{equation}\label{equ 1.4 UpperBoundofMDP}
\limsup_{N\rightarrow+\infty}\frac{N}{a_N^2}\log P\left(\theta^N\in C\right)\leq -\inf_{\pi\in C}\left(I_{ini}(\pi_0)+I_{dyn}(\pi)\right)
\end{equation}
for any closed set $C\subseteq \mathcal{D}([0, T_0], \mathcal{S})$ and
\begin{equation}\label{equ 1.5 LowerBoundofMDP}
\liminf_{N\rightarrow+\infty}\frac{N}{a_N^2}\log P\left(\theta^N\in O\right)\geq -\inf_{\pi\in O}\left(I_{ini}(\pi_0)+I_{dyn}(\pi)\right)
\end{equation}
for any open set $O\subseteq \mathcal{D}([0, T_0], \mathcal{S})$.
\end{theorem}
To make Theorem \ref{theorem 1.1 MDP} easy to catch, our next result gives alternative representation formulas of $I_{ini}$ and $I_{dyn}$. For any $f, g\in C([0, 1])$ and $0\leq t\leq T_0$, we define
\[
\langle f|g \rangle_t=\int_{[0,1]}\left(\int_0^1\lambda(x,y)\left(f(y)-f(x)\right)\left(g(y)-g(x)\right)dy\right)\mu_t(dx).
\]
Furthermore, for any $F, G\in C([0, T_0]\times [0, 1])$, we define
\[
\ll F, G\gg=\int_0^{T_0}\langle F_s|G_s \rangle_s ds.
\]
For $F_1, F_2\in C([0, T_0]\times [0, 1])$, we write $F_1\sim F_2$ when $\ll F_1-F_2, F_1-F_2\gg=0$. We use $\mathcal{H}$ to denote the completion of $C([0, T_0]\times[0, 1])/\sim$ under the inner product $\ll \cdot, \cdot\gg$, then we have the following result.
\begin{theorem}\label{theorem 1.3 alternative representation formula}
If $\nu\in \mathcal{S}$ makes $I_{ini}(\nu)<+\infty$, then there exists $g\in L^2([0, 1])$ that $\nu(dx)=g(x)\phi(x)dx$ and
\[
I_{ini}(\nu)=\nu(g)-\frac{1}{2}\int_0^1\phi(x)g^2(x)dx=\frac{\int_0^1\phi(x)g^2(x)dx}{2}.
\]
If $\pi\in \mathcal{D}([0, T_0], \mathcal{S})$ makes $I_{dyn}(\pi)<+\infty$, then there exists $F\in \mathcal{H}$ that
\[
\pi_{T_0}(G_{T_0})-\pi_0(G_0)-\int_0^{T_0}\pi_s\left((\partial_s+P_1-P_2)G_s\right)ds=\ll G, F\gg
\]
for any $G\in C^{1,0}([0, T_0]\times [0, 1])$ and $I_{dyn}(\pi)=\frac{\ll F, F\gg}{2}$.
\end{theorem}

Theorem \ref{theorem 1.3 alternative representation formula} is a routine result since $I_{ini}$ and $I_{dyn}$ are both defined as the supremum of a linear function minus a positive definite quadratic one. The proof of Theorem \ref{theorem 1.3 alternative representation formula} follows the same procedure as those in proofs of analogue results such as Lemma 5.1 of \cite{Kipnis1989} and Equation (2.2) of \cite{Xue2021}, where a crucial step is the utilization of Riesz representation theorem. Hence, we omit the proof of Theorem \ref{theorem 1.3 alternative representation formula} in this paper.

The proof of Theorem \ref{theorem 1.1 MDP} follows a routine strategy introduced in \cite{Kipnis1989}, where an exponential martingale plays the key role. A replacement lemma is crucial for the execution of the above strategy, which is the main difficulty we need to overcome in this paper. We prove this replacement lemma according to the large deviation principle of our model given in \cite{Xue2020}. For mathematical details, see Section \ref{section two}.

\section{Replacement lemma}\label{section two}
In this section we will prove the following replacement lemma.

\begin{lemma}\label{lemma 2.1 replacement}
Let $\mu$ be defined as in Equation \eqref{equation HydrodynamicLimit}, then for any $G\in C([0, T_0]\times[0, 1])$ and $\epsilon>0$,
\begin{equation}\label{equation supperexpoential in relacement lemma}
\limsup_{N\rightarrow+\infty}\frac{1}{a_N}\log P\left(\sup_{0\leq t\leq T_0}\left|\mu_t^N(G_t)-\mu_t(G_t)\right|\geq \epsilon\right)=-\infty.
\end{equation}
\end{lemma}
The large deviation principle of our model given in \cite{Xue2020} is crucial for the proof of Lemma \ref{lemma 2.1 replacement}, which we recall here. For any $\nu\in \mathcal{S}$, we define
\[
J_{ini}(\nu)=\sup_{f\in C([0, 1])}\left\{\nu(f)-\int_0^1\phi(x)\left(e^{f(x)}-1\right)dx\right\}.
\]
For any $\pi\in \mathcal{D}([0, T_0], \mathcal{S})$, we define
\[
J_{dyn}(\pi)=\sup_{G\in C^{1,0}([0, T_0]\times [0, 1])}\left\{\pi_{T_0}(G_{T_0})-\pi_0(G_0)-\int_0^{T_0}\pi_s\left((\partial_s+\mathcal{B})G_s\right)ds\right\},
\]
where
\[
\mathcal{B}f(x)=\int_0^1\lambda(x, y)\left(e^{f(y)-f(x)}-1\right)dy
\]
for any $f\in C([0, 1])$ and $x\in [0, 1]$. Then the following upper bound of large deviation principle is given in \cite{Xue2020}.

\begin{proposition}[{\cite[Theorem 2.6]{Xue2020}}]\label{proposition 2.2 LargeDeviation}
Let $\mu^N=\{\mu_t^N\}_{0\leq t\leq T_0}$, then
\[
\limsup_{N\rightarrow+\infty}\frac{1}{N}\log P\left(\mu^N\in C\right)\leq -\inf_{\pi\in C}(J_{ini}(\pi_0)+J_{dyn}(\pi))
\]
for any closed set $C\subseteq \mathcal{D}([0, T_0], \mathcal{S})$.
\end{proposition}
Note that although Reference \cite{Xue2020} adopts the assumption that $\lambda(x,y)=\lambda_1(x)\lambda_2(y)$ for some $\lambda_1, \lambda_2\in C([0, 1])$, this assumption is utilized in the proof of the lower bound of the large deviation principle. The upper bound does not rely on the this assumption.

To prove Lemma \ref{lemma 2.1 replacement}, we need the following two lemmas.

\begin{lemma}\label{lemma 2.3}
If $\pi\in \mathcal{D}([0, T_0], \mathcal{S})$ makes $J_{ini}(\pi_0)+J_{dyn}(\pi)=0$, then $\pi=\mu$.
\end{lemma}

\begin{lemma}\label{lemma 2.4}
For any $0<C<+\infty$,
\[
\mathfrak{A}_C:=\left\{\pi\in \mathcal{D}([0, T_0], \mathcal{S}):~J_{ini}(\pi_0)+J_{dyn}(\pi)\leq C\text{~and~}\pi_t\in \mathcal{A}\text{~for all~}0\leq t\leq T_0\right\}
\]
is compact.
\end{lemma}
We first utilize Lemmas \ref{lemma 2.3} and \ref{lemma 2.4} to prove Lemma \ref{lemma 2.1 replacement}.

\proof[Proof of Lemma \ref{lemma 2.1 replacement}]

For any $\epsilon>0$ and given $G\in C^{1,0}([0, T_0]\times[0, 1])$, we define $D_{\epsilon, G}$ as
\[
D_{\epsilon, G}=\left\{\pi\in \mathcal{D}([0, T_0], \mathcal{S}):~\sup_{0\leq t\leq T_0}\left|\pi_t(G_t)-\mu_t(G_t)\right|\geq \epsilon\text{~and~}\pi_t\in \mathcal{A}\text{~for all~}0\leq t\leq T_0\right\}.
\]
Since $\mu^N_t\in \mathcal{A}$ for all $0\leq t\leq T_0$ and $\frac{N^2}{a_N^2}\rightarrow+\infty$, by Proposition \ref{proposition 2.2 LargeDeviation}, we only need to show that
\[
\inf_{\pi\in D_{\epsilon, G}}\left(J_{ini}(\pi_0)+J_{dyn}(\pi)\right)>0
\]
to prove Lemma \ref{lemma 2.1 replacement}. If $\inf_{\pi\in D_{\epsilon, G}}\left(J_{ini}(\pi_0)+J_{dyn}(\pi)\right)=0$, then there exists a sequence $\{\pi^n\}_{n\geq 1}$ in $D_{\epsilon, G}\cap \mathfrak{A}_1$ that
\begin{equation}\label{equ 2.4}
\lim_{n\rightarrow+\infty}\left(J_{ini}(\pi^n_0)+J_{dyn}(\pi^n)\right)=0.
\end{equation}
By Lemma \ref{lemma 2.4}, $\mathfrak{A}_1$ is compact. Hence, there exists $\hat{\pi}\in \mathfrak{A}_1$ that a subsequence $\{\pi^{n_k}\}_{k\geq 1}$ of $\{\pi^n\}_{n\geq 1}$ satisfies that $\lim_{k\rightarrow+\infty}\pi^{n_k}=\hat{\pi}$. Since $J_{ini}$ and $J_{dyn}$ are both defined as supremums of continuous functions, it is easy to check that $J_{ini}(\pi_0)+J_{dyn}(\pi)$ is lower semi-continuous of $\pi$. Then, by Equation \eqref{equ 2.4},
\[
J_{ini}(\hat{\pi}_0)+J_{dyn}(\hat{\pi})=0
\]
and consequently $\hat{\pi}=\mu$ according to Lemma \ref{lemma 2.3}. However, since $D_{\epsilon, G}$ is closed, $\hat{\mu}\in D_{\epsilon, G}$ and hence
\[
\sup_{0\leq t\leq T_0}\left|\hat{\pi}_t(G_t)-\mu_t(G_t)\right|\geq \epsilon,
\]
which is contradict with $\hat{\pi}=\mu$.

\qed

At last we prove Lemmas \ref{lemma 2.3} and \ref{lemma 2.4}.

\proof[Proof of Lemma \ref{lemma 2.3}]

Since $J_{ini}(\pi_0)\geq \pi_0(0)-\int_0^1\phi(x)\left(e^0-1\right)dx=0$ and
\[
J_{dyn}(\pi)\geq\pi_{T_0}(0)-\pi_0(0)-\int_0^{T_0}\pi_s\left((\partial_s+\mathcal{B})0\right)ds=0,
\]
$J_{ini}(\pi_0)+J_{dyn}(\pi)=0$ implies that $J_{ini}(\pi_0)=J_{dyn}(\pi)=0$. Then, for any $f\in C([0, 1])$,
\[
K_f(c):=\pi_0(cf)-\int_0^1\phi(x)\left(e^{cf(x)}-1\right)dx
\]
gets maximum $0$ at $c=0$ and hence $\frac{d}{dc}K_f(c)\Big|_{c=0}=0$, which implies that
\[
\pi_0(f)=\int_0^1\phi(x)f(x)dx
\]
for any $f\in C([0,1])$ and hence $\pi_0(dx)=\phi(x)dx$. Similarly, for any $h\in C^1([0,T_0])$ and $f\in C([0, 1])$, let $G^{h,f}(t,x)=h(t)f(x)$, then
\[
\Gamma_{h,f}(c):=\pi_{T_0}(cG^{h,f}_{T_0})-\pi_0(cG^{h,f}_0)-\int_0^{T_0}\pi_s\left((\partial_s+\mathcal{B})(cG^{h,f}_s)\right)ds
\]
gets maximum at $c=0$ and hence $\frac{d}{dc}\Gamma_{h,f}(c)\Big|_{c=0}=0$, which implies that
\[
h_{T_0}\pi_{T_0}(f)-h_0\pi_0(f)-\int_0^{T_0}h^\prime_s\pi_s(f)ds=\int_0^{T_0}h_s\pi_s\left((P_1-P_2)f\right)ds
\]
for any $f\in C([0, 1]), h\in C^1([0, T_0])$ and hence $\{\pi_t(f)\}_{0\leq t\leq T_0}$ is differentiable with
\[
\frac{d}{dt}\pi_t(f)=\pi_t\left((P_1-P_2)f\right)
\]
for any $f\in C([0, 1])$. Consequently, $\pi=\mu$.

\qed

\proof[Proof of Lemma \ref{lemma 2.4}]
By Arzel\`{a}-Ascoli Theorem, we only need to show that for any nonnegative $f\in C[0, 1]$,
\[
\left\{\pi_t(f):~0\leq t\leq T_0\right\}_{\pi\in \mathfrak{A}_C}
\]
are uniformly bounded and equicontinuous. Let $\vec{1}$ be the function that $\vec{1}(x)=1$ for all $x\in [0, 1]$, then for any $\pi\in\mathfrak{A}_C$,
\begin{equation}\label{equ 2.2}
\pi_0(\vec{1})\leq C+\int_0^1\phi(x)\left(e^{\vec{1}(x)}-1\right)dx=C+(e-1)\int_0^1\phi(x)dx.
\end{equation}
For given $0<t<T_0$ and sufficiently large $n$, let $\Lambda^t_n$ be the function from $[0, T_0]$ to $\mathbb{R}$ that $\Lambda^t_n(s)=0$ when $s\leq t$ or $s\geq t+\frac{1}{n}$ and $\Lambda^t_n(s)=-n$ when $t<s<t+\frac{1}{n}$. Since $C([0, T_0])$ is dense in $L^1([0, T_0])$, let $\{\tilde{\Lambda}^t_{n,m}\}_{m\geq 1}$ be a sequence in $C([0, T_0])$ that $\tilde{\Lambda}^t_{n,m}$ converges in $L^1$ to $\Lambda^t_n$ as $m\rightarrow+\infty$. Then we define $h^t_n\in C([0,T_0]), \tilde{h}^t_{n,m}\in C^1([0, T_0])$ that
\[
h^t_n(s)=1+\int_0^s\Lambda^t_n(u)du,~\tilde{h}^t_{n,m}(s)=1+\int_0^s\tilde{\Lambda}^t_{n,m}(u)du
\]
for all $s\in [0, 1]$. As a result, $\tilde{h}^t_{n,m}$ converges to $h^t_n$ uniformly in $[0, 1]$ as $m\rightarrow+\infty$. Let $\tilde{G}^t_{n,m}\in C^{1,0}([0, T_0]\times[0, 1])$ that $\tilde{G}^t_{n,m}(s, x)=\tilde{h}^t_{n,m}(s)\vec{1}(x)$ for any $0\leq s\leq T_0, 0\leq x\leq 1$, then for $\pi\in \mathfrak{A}_C$,
\[
\pi_{T_0}\left(\tilde{G}^t_{n,m,T_0}\right)-\pi_{0}\left(\tilde{G}^t_{n,m,0}\right)\leq C+\int_0^{T_0}\pi_s\left((\partial_s+\mathcal{B})\tilde{G}^t_{n,m,s}\right)ds.
\]
Let $m\rightarrow+\infty$, we have
\[
-\pi_0(\vec{1})\leq C-n\int_t^{t+\frac{1}{n}}\pi_s(\vec{1})ds.
\]
Since $\pi$ is right-continuous, let $n\rightarrow\infty$, we have
\begin{equation}\label{equ 2.1point5}
\pi_t(\vec{1})\leq \pi_0(\vec{1})+C.
\end{equation}
Then, by Equation \eqref{equ 2.2},
\begin{align}\label{equ 2.3}
\pi_t(f)&\leq \left(\max_{0\leq x\leq 1}f(x)\right)\pi_t(\vec{1}) \notag\\
&\leq \left(\max_{0\leq x\leq 1}f(x)\right)\left(2C+(e-1)\int_0^1\phi(x)dx\right)
\end{align}
for any $0\leq t\leq T_0$ and hence $\left\{\pi_t(f):~0\leq t\leq T_0\right\}_{\pi\in \mathfrak{A}_C}$ are uniformly bounded. For $s<t<T_0$ and sufficiently large $n$, we define $\hat{\Lambda}^{t,s}_n$ as the function from $[0, T_0]$ to $\mathbb{R}$ that
\[
\hat{\Lambda}^{t,s}_n(u)=
\begin{cases}
0 &\text{~if~}u\leq s, s+\frac{1}{n}<u\leq t \text{~or~}u>t+\frac{1}{n},\\
n &\text{~if~}s<u\leq s+\frac{1}{n},\\
-n &\text{~if~}t<u\leq t+\frac{1}{n}.
\end{cases}
\]
Then, via replacing $\Lambda^t_n$ by $\hat{\Lambda}^{t,s}_n$ and $\vec{1}$ by $g\in C([0, 1])$ in the analysis leading to Equation \eqref{equ 2.1point5}, we have
\[
\pi_t(g)\leq \pi_s(g)+C+\int_s^t\pi_u\left(\mathcal{B}g\right)du
\]
for any $g\in C([0, 1])$. For any given $M>0$, let $g=Mf$, then we have
\[
\pi_t(f)\leq \pi_s(f)+\frac{C}{M}+\int_s^t \frac{1}{M}\pi_u\left(\mathcal{B}(Mf)\right)du.
\]
By Equation \eqref{equ 2.3}, for any $\pi\in \mathfrak{A}_C$ and $u\in (s, t)$,
\[
\pi_u\left(\mathcal{B}(Mf)\right)\leq \left(\max_{0\leq x\leq 1}|\mathcal{B}(Mf)(x)|\right)\left(2C+(e-1)\int_0^1\phi(x)dx\right).
\]
Hence, for any $\epsilon>0$, we can first choose $M$ sufficiently large that $\frac{C}{M}<\frac{\epsilon}{2}$ and then there exists $\delta_1>0$ only depending on $M$ and $f$ that
\[
\int_s^t \frac{1}{M}\pi_u\left(\mathcal{B}(Mf)\right)du<\frac{\epsilon}{2}
\]
and hence
\[
\pi_t(f)\leq \pi_s(f)+\epsilon
\]
for any $t-s\leq \delta_1$ and $\pi\in \mathfrak{A}_C$. Let $g=-Mf$, then it is proved similarly that there exists $\delta_2>0$ only depending on $\epsilon$ and $f$ that
\[
\pi_t(f)\geq \pi_s(f)-\epsilon
\]
for any $t-s\leq \delta_2$ and $\pi\in \mathfrak{A}_C$. As a result, $\{\pi_t(f):~0\leq t\leq T_0\}_{\pi\in \mathfrak{A}_C}$ are equicontinuous and hence the proof is complete.

\qed

\section{The proof of Equation \eqref{equ 1.4 UpperBoundofMDP}}\label{Section 3}
In this section we give the proof of Equation \eqref{equ 1.4 UpperBoundofMDP}. With Lemma \ref{lemma 2.1 replacement}, the proof of our main result follows a routine procedure introduced in \cite{Kipnis1989}, which has also been utilized in References \cite{Gao2003}, \cite{Xue2021} and so on to prove MDPs of models such as exclusion processes, density-dependent Markov chains and so on. Hence, in this paper we only give a outline of the proof without repeating too many similar details with those in above references.

For later use, for a given positive sequence $\{c_N\}_{N\geq 1}$ that $\lim_{N\rightarrow+\infty}c_N=+\infty$ and a sequence of random variables $\{Y_N\}_{N\geq 1}$, we write $Y_N$ as $o_{\exp}(c_N)$ when
\[
\lim_{N\rightarrow+\infty}\frac{1}{c_N}\log P(|Y_N|\geq \epsilon)=-\infty
\]
for any $\epsilon>0$ and write $Y_N$ as $O_{\exp}(c_N)$ when
\[
\limsup_{N\rightarrow+\infty}\frac{1}{c_N}\log P(|Y_N|\geq \epsilon)<0
\]
for any $\epsilon>0$. Now we first prove Equation \eqref{equ 1.4 UpperBoundofMDP} for compact $K\subseteq \mathcal{D}([0, T_0], \mathcal{S})$.

\proof [Proof of Equation \eqref{equ 1.4 UpperBoundofMDP} for compact sets]

For each $N\geq 1$ and any $G\in C^{1,1}([0, T_0]\times [0, 1])$, we define $H^N_G(t, X_t^N)$ as
\[
H_G^N(t, X_t^N)=\exp\left\{\frac{a_N^2}{N}\theta_t^N(G_t)\right\}
=\exp\left\{\frac{a_N}{N}\sum_{i=1}^N\left(X_t^N(i)-EX_t^N(i)\right)G_t(\frac{i}{N})\right\}
\]
and define $\Gamma_t^N(G)$ as
\[
\Gamma_t^N(G)=\frac{H_G^N(t, X_t^N)}{H_G^N(0, X_0^N)}\exp\left\{-\int_0^t\frac{\left(\partial_s+\mathcal{L}_N\right)H_G^N(s, X_s^N)}{H_G^N(s, X_s^N)}ds\right\}.
\]
Then it is easy to check that $\{\Gamma_t^N(G)\}_{0\leq t\leq T_0}$ is a martingale with mean $1$ by It\^{o}'s formula. Therefore, for any $t\geq 0$ and $f\in C([0, 1])$,
\begin{equation}\label{equ 3.1}
Ee^{\frac{a_N}{N}\sum_{i=1}^N\left(X_0^N(i)-EX_0^N(i)\right)f(\frac{i}{N})}
=E\left(e^{\frac{a_N}{N}\sum_{i=1}^N\left(X_0^N(i)-EX_0^N(i)\right)f(\frac{i}{N})}\Gamma_t^N(G)\right).
\end{equation}
According to our assumption of $X_0^N$ and the fact that $\lim_{N\rightarrow+\infty}\frac{a_N}{N}=0$, it is easy to check that
\begin{equation}\label{equ 3minus1}
\lim_{N\rightarrow+\infty}Ee^{\frac{a_N}{N}\sum_{i=1}^N\left(X_t^N(i)-EX_0^N(i)\right)f(\frac{i}{N})}=\frac{1}{2}\int_0^1\phi(x)f^2(x)dx.
\end{equation}
according to Taylor's expansion formula up to the second order.
For later use, for each $N\geq 1$, we define $P_1^Nf(x)=\frac{1}{N}\sum_{j=1}^N\lambda(x, \frac{j}{N})f(\frac{j}{N})$, $P_2^Nf(x)=f(x)\frac{1}{N}\sum_{j=1}^N\lambda(x, \frac{j}{N})$,
$\mathcal{K}f(x)=\int_0^1\lambda(x, y)(f(y)-f(x))^2dy$ and $\mathcal{K}^Nf(x)=\frac{1}{N}\sum_{j=1}^N\lambda(x, \frac{j}{N})(f(\frac{j}{N})-f(x))^2$ for any $f\in C([0, 1]), x\in [0, 1]$.
According to the generator $\mathcal{L}_N$ of $\{X_t^N\}_{t\geq 0}$,
\[
\frac{d}{dt}EX_t^N(i)=-EX_t^N(i)\sum_{j=1}^N\frac{\lambda(\frac{i}{N},\frac{j}{N})}{N}+\sum_{j=1}^N\frac{\lambda(\frac{j}{N}, \frac{i}{N})}{N}EX_t^N(j)
\]
while
\[
\mathcal{L}_NH_G^N(t, X_t^N)=\sum_{i=1}^N\frac{\lambda(\frac{i}{N}, \frac{j}{N})X_t^N(i)}{N}\sum_{j=1}^NH_G^N(t, X_t^N)\left(e^{\frac{a_N}{N}\left(G_t(\frac{j}{N})-G_t(\frac{i}{N})\right)}-1\right).
\]
Then, by the fact that $\frac{a_N}{N}\rightarrow 0$ and Taylor's expansion formula up to the second order, it is not difficult to show that
\begin{equation}\label{equ 3.0}
\Gamma_{T_0}^N(G)=\exp\left\{\frac{a_N^2}{N}\left(l(\theta^N, G)+\epsilon^N\right)\right\},
\end{equation}
where
\begin{align*}
l(\pi, G)=&\pi_{T_0}(G_{T_0})-\pi_0(G_0)-\int_0^{T_0}\pi_s\left((\partial_s+P_1-P_2)G_s\right)ds  \\ &-\frac{1}{2}\int_0^{T_0}\left(\int_{[0,1]}\left(\int_0^1\lambda(x,y)\left(G_s(y)-G_s(x)\right)^2dy\right)\mu_s(dx)\right)ds \notag
\end{align*}
for any $\pi\in \mathcal{D}([0, T_0], \mathcal{S})$ and
\[
\epsilon^N=\int_0^{T_0}\left(\epsilon_{1, t}^N+\epsilon_{2, t}^N+\epsilon_{3,t}^N+\epsilon_{4, t}^N\right)dt,
\]
where $\epsilon_{1,t}^N$ is the third order Lagrange's remainder of the Taylor's formula that
\[
|\epsilon_{1, t}^N|\leq C_1\frac{a_N}{N}\frac{1}{N}\sum_{i=1}^N\left(X_t^N(i)+EX_t^N(i)\right)
\]
with constant $C_1<+\infty$ independent of $t$ and $N$,
\[
\epsilon_{2,t}^N=\mu^N_t(\mathcal{K}G_t)-\mu_t(\mathcal{K}G_t), \text{~}\epsilon_{3,t}^N=\mu^N_t(\mathcal{K}^NG_t)-\mu^N_t(\mathcal{K}G_t)
\]
and
\[
\epsilon_{4,t}^N=\theta_t^N\left((P_1^N-P_2^N)G_t-(P_1-P_2)G_t\right).
\]
According to the fact that $\sum_{i=1}^NX_t^N(i)\equiv \sum_{i=1}^NX_0^N(i)$ and our assumption of $X_0^N$, it is easy to check that $\sup_{t\leq T_0}|\epsilon_{1, t}^N|=O_{\exp}\left(\frac{N^2}{a_N}\right)=o_{\exp}\left(a_N\right)$ by Markov's inequality. Since $\lambda\in C^{1,1}([0, 1]\times[0,1]), G_t\in C^1([0, 1])$, it is easy to check that
\[
|\epsilon_{4,t}^N|\leq \frac{C_2}{Na_N}\sum_{i=1}^N(X_0^N(i)+EX_0^N(i))
\]
for some $C_2<+\infty$ independent of $t$ and $N$ according to Lagrange's mean value theorem. Then, we similarly have $\sup_{t\leq T_0}|\epsilon_{4, t}^N|=O_{\exp}\left(Na_N\right)=o_{\exp}\left(a_N\right)$ according to Markov's inequality. According to a similar analysis,
\[
|\epsilon_{3,t}^N|\leq \frac{C_3}{N^2}\sum_{i=1}^NX_0^N(i)
\]
for some $C_3<+\infty$ independent of $t, N$ and hence $\sup_{t\leq T_0}|\epsilon_{3, t}^N|=O_{\exp}\left(N^2\right)=o_{\exp}\left(a_N\right)$ according to Markov's inequality. By Lemma \ref{lemma 2.1 replacement}, $\sup_{t\leq T_0}|\epsilon_{2, t}^N|=o_{\exp}\left(a_N\right)$. In conclusion,
\begin{equation*}
\epsilon^N=o_{\exp}\left(a_N\right)=o_{\exp}\left(\frac{a_N^2}{N}\right).
\end{equation*}
As a result, for any $\epsilon>0$ and compact $K\subseteq \mathcal{D}\big([0, T_0], \mathcal{S}\big)$,
\begin{equation}\label{equ 3.3}
\limsup_{N\rightarrow+\infty}\frac{N}{a_N^2}\log P\left(\theta^N\in K, |\epsilon^N|\leq \epsilon\right)
=\limsup_{N\rightarrow+\infty}\frac{N}{a_N^2}\log P\left(\theta^N\in K\right).
\end{equation}
By Equation \eqref{equ 3.0}, $\Gamma_{T_0}^N(G)\geq \exp\left\{\frac{a_N^2}{N}\left(l(\theta^N, G)-\epsilon\right)\right\}$ when $|\epsilon^N|\leq \epsilon$. Therefore, by Equation \eqref{equ 3.1},
\begin{align*}
&Ee^{\frac{a_N}{N}\sum_{i=1}^N\left(X_0^N(i)-EX_0^N(i)\right)f(\frac{i}{N})}\\
&\geq E\left(e^{\frac{a_N}{N}\sum_{i=1}^N\left(X_0^N(i)-EX_0^N(i)\right)f(\frac{i}{N})}\Gamma_{T_0}^N(G)1_{\{\theta^N\in K, |\epsilon^N|\leq \epsilon\}}\right)\\
&\geq \exp\left\{\frac{a_N^2}{N}\inf_{\pi\in K}\left\{\pi_0(f)+l(\pi, G)-\epsilon\right\}\right\}P\left(\theta^N\in K, |\epsilon^N|\leq \epsilon\right).
\end{align*}
Then, according to Equations \eqref{equ 3minus1} and \eqref{equ 3.3},
\begin{align*}
&\limsup_{N\rightarrow+\infty}\frac{N}{a_N^2}\log P\left(\theta^N\in K\right)
=\limsup_{N\rightarrow+\infty}\frac{N}{a_N^2}\log P\left(\theta^N\in K, |\epsilon^N|\leq \epsilon\right)\\
&\leq -\inf_{\pi\in K}\left\{\pi_0(f)+l(\pi, G)\right\}+\frac{1}{2}\int_0^1\phi(x)f^2(x)dx+\epsilon\\
&=-\inf_{\pi\in K}\left\{\pi_0(f)-\frac{1}{2}\int_0^1\phi(x)f^2(x)dx+l(\pi, G)\right\}+\epsilon.
\end{align*}
Since $f, G, \epsilon$ are arbitrary,
\begin{align}\label{equ 3.4}
&\limsup_{N\rightarrow+\infty}\frac{N}{a_N^2}\log P\left(\theta^N\in K\right) \notag\\
&\leq -\sup_{f\in C([0, 1]), \atop G\in C^{1,1}([0, T_0]\times[0, 1])}\inf_{\pi\in K}\left\{\pi_0(f)-\frac{1}{2}\int_0^1\phi(x)f^2(x)dx+l(\pi, G)\right\}.
\end{align}
Since $\pi_0(f)-\frac{1}{2}\int_0^1\phi(x)f^2(x)dx+l(\pi, G)$ is concave with $(f, G)$ and convex with $\pi$, according to the minimax theorem given in
\cite{Sion1958},
\begin{align*}
&\sup_{f\in C([0, 1]), \atop G\in C^{1,1}([0, T_0]\times[0, 1])}\inf_{\pi\in K}\left\{\pi_0(f)-\frac{1}{2}\int_0^1\phi(x)f^2(x)dx+l(\pi, G)\right\}\\
&=\inf_{\pi\in K}\sup_{f\in C([0, 1]), \atop G\in C^{1,1}([0, T_0]\times[0, 1])}\left\{\pi_0(f)-\frac{1}{2}\int_0^1\phi(x)f^2(x)dx+l(\pi, G)\right\}\\
&=\inf_{\pi \in K}\left(\sup_{f\in C([0, 1])}\left\{\pi_0(f)-\frac{1}{2}\int_0^1\phi(x)f^2(x)dx\right\}+\sup_{G\in C^{1,1}([0, T_0]\times[0, 1])}l(\pi, G)\right)\\
&=\inf_{\pi \in K}\left(I_{ini}(\pi_0)+\sup_{G\in C^{1,1}([0, T_0]\times[0, 1])}l(\pi, G)\right).
\end{align*}
Since $C^{1,1}\left([0, T_0]\times[0,1]\right)$ is dense in $C^{1,0}\left([0, T_0]\times[0, 1]\right)$,
\[
\sup_{G\in C^{1,1}([0, T_0]\times[0, 1])}l(\pi, G)=\sup_{G\in C^{1,0}([0, T_0]\times[0, 1])}l(\pi, G)=I_{dyn}(\pi)
\]
and hence Equation \eqref{equ 1.4 UpperBoundofMDP} holds for all compact $K\subseteq \mathcal{D}\left([0, T_0], \mathcal{S}\right)$ according to Equation \eqref{equ 3.4}.

\qed

To prove Equation \eqref{equ 1.4 UpperBoundofMDP} for all closed sets, we need the following two lemmas as preliminaries.

\begin{lemma}\label{lemma 3.1 independent and Poisson}
Under our assumption of $X_0^N$, $X_t^N(1), X_t^N(2),\ldots, X_t^N(N)$ are independent for any $t\geq 0$ and $X_t^N(i)$ follows Poisson distribution with mean $EX_t^N(i)$ for all $1\leq i\leq N$.
\end{lemma}

\begin{lemma}\label{lemma 3.2 control of integration}
For any $f\in C([0, 1])$ and $\epsilon>0$,
\[
\limsup_{M\rightarrow+\infty}\limsup_{N\rightarrow+\infty}\frac{N}{a_N^2}\log P\left(\sup_{0\leq t\leq T_0}\left|\int_0^t\theta^N_s(f)ds\right|>M\right)=-\infty
\]
and
\[
\limsup_{\delta\rightarrow0}\limsup_{N\rightarrow+\infty}\frac{N}{a_N^2}\log P\left(\sup_{|t-s|\leq\delta\atop 0\leq s<t\leq T_0}\left|\int_s^t\theta_u^N(f)du\right|>\epsilon\right)=-\infty.
\]
\end{lemma}

The proof of Lemma \ref{lemma 3.1 independent and Poisson} is given in Appendix \ref{subsection A.1}. With Lemma \ref{lemma 3.1 independent and Poisson}, we have
\begin{equation}\label{equ 3.6}
E\left(\exp\{\frac{a^2_N}{N}\theta_t^N(f)\}\right)=\exp\left\{\sum_{i=1}^NEX_t^N(i)\left(e^{\frac{a_N}{N}f\left(\frac{i}{N}\right)}
-\frac{a_N}{N}f\left(\frac{i}{N}\right)-1\right)\right\}.
\end{equation}
We have shown in \cite{Xue2020} that there exists $C_5<+\infty$ independent of $N$ that
\begin{equation}\label{equ 3.7}
\sup_{N\geq 1, 1\leq i\leq N, 0\leq t\leq T_0}EX_t^N(i)\leq C_5.
\end{equation}
With Equations \eqref{equ 3.6} and \eqref{equ 3.7}, the proof of Lemma \ref{lemma 3.2 control of integration} follows the same procedure as that introduced in the proof of Lemma 2.2 of \cite{Gao2003}, where a crucial step is the utilization of Garsia-Rademich-Rumsey Lemma. Hence we omit the details of the proof of Lemma \ref{lemma 3.2 control of integration} here.

At last, we give the proof of Equation \eqref{equ 1.4 UpperBoundofMDP} for all closed sets.

\proof[Proof of Equation \eqref{equ 1.4 UpperBoundofMDP}]

Since we have proved Equation \eqref{equ 1.4 UpperBoundofMDP} for all compact sets, we only need to show that $\{\theta^N\}_{N\geq 1}$ are exponentially tight to complete this proof. By the criteria given in \cite{Puhalskii1994}, we only need to show that
\begin{equation}\label{equ 3.8}
\limsup_{M\rightarrow+\infty}\limsup_{N\rightarrow+\infty}\frac{N}{a_N^2}\log P\left(\sup_{0\leq t\leq T_0}\left|\theta_t^N(f)\right|>M\right)=-\infty
\end{equation}
and
\begin{equation}\label{equ 3.9}
\limsup_{\delta\rightarrow 0}\limsup_{N\rightarrow+\infty}\frac{N}{a_N^2}\log \sup_{\tau\in \Upsilon}P\left(\sup_{0<t\leq \delta}\left|\theta^N_{\tau+t}(f)-\theta_\tau^N(f)\right|>\epsilon\right)=-\infty
\end{equation}
for any $f\in C^1([0,1])$ and $\epsilon>0$, where $\Upsilon$ is the set of all stopping times bounded by $T_0$ from above. With Lemma \ref{lemma 3.2 control of integration}, proofs of Equations \eqref{equ 3.8} and \eqref{equ 3.9} follows same procedures as those in proofs of Equations (3.3) and (3.4) of \cite{Gao2003} respectively, where a crucial step is the utilization of Doob's inequality on the exponential martingale $\{\Gamma_t^N(f)\}_{0\leq t\leq T_0}$. Consequently, the proof is complete.

\qed

\section{The proof of Equation \eqref{equ 1.5 LowerBoundofMDP}}\label{section 4}

In this section we prove Equation \eqref{equ 1.5 LowerBoundofMDP}. As we have introduced, our proof follows the strategy introduced in \cite{Kipnis1989}, where a crucial step is to derive the law of large numbers of $\theta^N$ under the transformed probability measure with $\Gamma^N_{T_0}(G)$ introduced in Section \ref{Section 3} as the R-N derivative with respect to the original measure of $\{X_t^N\}_{t\geq 0}$. For later use, we first introduce some notations and definitions. For any $f\in C([0,1])$ and sufficiently large $N$, we denote by $P^N_f$ the probability measure of our process $\{X_t^N\}_{t\geq 0}$ under the initial condition that $\{X_0^N(i)\}_{1\leq i\leq N}$ are independent and $X_0^N(i)$ follows Poisson distribution with mean $\phi\left(\frac{i}{N}\right)+\frac{a_N}{N}f\left(\frac{i}{N}\right)$. For any $G\in C^{1,1}\left([0, T_0]\times [0, 1]\right)$, we define $\hat{P}^N_{f, G}$ as the probability measure that
\[
\frac{d\hat{P}^N_{f, G}}{dP^N_f}=\Gamma_{T_0}^N(G).
\]
Then the following lemma is crucial for us to prove Equation \eqref{equ 1.5 LowerBoundofMDP}, which gives the law of large numbers of $\theta^N$ under the transformed measure $\hat{P}^N_{f, G}$.
\begin{lemma}\label{lemma 4.1 LLNunderTransformedMeasure}
For given $G\in C^{1,1}\left([0, T_0]\times [0, 1]\right)$, $\theta^N$ converges in $\hat{P}^N_{f, G}$-probability to $\vartheta^{f,G}$ as $N\rightarrow+\infty$, where $\vartheta^{f, G}$ is the unique element in $\mathcal{D}\left([0, T_0], \mathcal{S}\right)$ that
\begin{equation}\label{equ 4.1 measureValuedODE}
\begin{cases}
\frac{d}{dt}\vartheta_t^{f,G}(h)=\vartheta_t^{f, G}\left((P_1-P_2)h\right)+\langle G_t | h\rangle_t \text{~for all~}0\leq t\leq T_0\text{~and~}h\in C([0, 1]),\\
\vartheta_0^{f, G}(dx)=f(x)dx.
\end{cases}
\end{equation}
\end{lemma}
\begin{remark}\label{Remark 4.1}
Lemma \ref{lemma 4.1 LLNunderTransformedMeasure} is a routine auxiliary result for the proof of the lower bound of the MDP. Analogues of Lemma \ref{lemma 4.1 LLNunderTransformedMeasure} such as Theorems 3.1 of \cite{Kipnis1989}, 4.1 of \cite{Gao2003} and Lemma 4.2 of \cite{Xue2021} have been given in literatures to prove LDPs or MDPs for models such as exclusion processes and density dependent Markov chains. With Lemma \ref{lemma 4.1 LLNunderTransformedMeasure}, roughly speaking, we can estimate $P(\theta^N=d\pi)$ for some $\pi\in \mathcal{D}\left([0, T_0], \mathcal{S}\right)$ as following. Choose $f, G$ to make $\vartheta^{f,G}=\pi$, then
\begin{align*}
P(\theta^N=d\pi)=E_{\hat{P}^N_{f, G}}\left(\frac{dP}{dP^N_f}\left(\Gamma_{T_0}^N(G)\right)^{-1}1_{\{\theta^N=d\pi\}}\right).
\end{align*}
Lemma \ref{lemma 4.1 LLNunderTransformedMeasure} implies that $\hat{P}^N_{f, G}\left(\theta^N=d\pi\right)=1+o(1)$ and hence our MDP holds when we can show that
\[
\frac{dP}{dP^N_f}\left(\Gamma_{T_0}^N(G)\right)^{-1}\Bigg|_{\theta^N=\pi}=\exp\left\{-\frac{a_N^2}{N}\left(I_{ini}(\pi_0)+I_{dyn}(\pi)+o(1)\right)\right\},
\]
which can be obtained according to Theorem \ref{theorem 1.3 alternative representation formula}. The rigorous statement of the above intuitive analysis is given at the end of this section.
\end{remark}

The following lemma is a preliminary for us to prove Lemma \ref{lemma 4.1 LLNunderTransformedMeasure}.
\begin{lemma}\label{lemma 4.2}
For given $f, h\in C([0, 1])$ and $G\in C^{1,1}([0, T_0]\times[0, 1])$,
\begin{equation}\label{equ 4.2}
\sum_{0\leq t\leq T_0}\left(\theta^N_t(h)-\theta^N_{t-}(h)\right)^2=o_{\exp}(N)
\end{equation}
under both $P_f^N$ and $\hat{P}_{f,G}^N$.
\end{lemma}
The proof of Lemma \ref{lemma 4.2} is given in Appendix \ref{subsection A.2}. For the proof of Lemma \ref{lemma 4.1 LLNunderTransformedMeasure}, we introduce some notations and definitions. For a sequence of random variables $\{Y_N\}_{N\geq 1}$, we write $Y_N$ as $o_p(1)$ when $\lim_{N\rightarrow+\infty}Y_N=0$ in probability. For any $h\in C^1([0, 1])$ and $G\in C^{1,1}([0, T_0]\times [0, 1])$, we define
\begin{align*}
\mathcal{M}_t(\theta^N(h))&=\theta_t^N(h)-\theta_0^N(h)-\int_0^t\left(\partial_s+\mathcal{L}_N\right)\theta_s^N(h)ds\\
&=\theta_t^N(h)-\theta_0^N(h)-\int_0^t\theta_s^N\left((P_1^N-P_2^N)h\right)ds
\end{align*}
and
\[
\mathcal{M}_t(H_G^N)=H_G^N(t, X_t^N)-H_G^N(0, X_0^N)-\int_0^t\left(\partial_s+\mathcal{L}_N\right)H_G^N(s, X_s^N)ds,
\]
where $H_G^N$ is defined as in Section \ref{Section 3}. According to basic properties of Markov processes, $\{\mathcal{M}_t(\theta^N(h))\}_{t\geq 0}$ and $\{\mathcal{M}_t(H_G^N)\}_{t\geq 0}$ are both martingales. In this paper, for two local martingales $\{\mathcal{M}^1_t\}_{t\geq 0}, \{\mathcal{M}^2\}_{t\geq 0}$, we use
$\{\langle \mathcal{M}^1, \mathcal{M}^2\rangle_t\}_{t\geq 0}$ to denote the predictable quadratic-covariation process which is continuous and use $\{[\mathcal{M}^1, \mathcal{M}^2]_t\}_{t\geq 0}$ to denote the optional quadratic-covariation process that
\[
\lim_{\sup_i(t_{i+1}-t_i)\rightarrow 0}\sum_{i}\left(\mathcal{M}^1_{t_{i+1}}-\mathcal{M}^1_{t_i}\right)\left(\mathcal{M}^2_{t_{i+1}}-\mathcal{M}^2_{t_i}\right)=[\mathcal{M}^1, \mathcal{M}^2]_t
\]
in probability, where the limit is over all partitions $\{t_i\}$ of $[0, t]$. Then, according to basic properties of Markov processes and direct calculations,
\begin{align}\label{equ 4.3}
&d\langle \mathcal{M}(H_G^N), \mathcal{M}(\theta^N(h))\rangle_t \notag\\
&=\left(-\theta_t^N(h)\mathcal{L}_NH_G^N(t, X_t^N)-H_G^N(t, X_t^N)\mathcal{L}_N\theta_t^N(h)+\mathcal{L}_N\left(\theta^N_t(h)H_G^N(t, X_t^N)\right)\right)dt\\
&=\sum_{i=1}^N\sum_{j=1}^N\frac{\lambda\left(\frac{i}{N}, \frac{j}{N}\right)X_t^N(i)}{N}H_G^N(t, X_t^N)\left(e^{\frac{a_N}{N}\left(G_t\left(\frac{j}{N}\right)-G_t\left(\frac{i}{N}\right)\right)}-1\right)\frac{h\left(\frac{j}{N}\right)
-h\left(\frac{i}{N}\right)}{a_N}dt\notag.
\end{align}

Now we prove Lemma \ref{lemma 4.1 LLNunderTransformedMeasure}. Our proof follows the strategy introduced in the proof of Lemma 4.2 of \cite{Xue2021}, where a crucial step is the utilization of a generalized version of Girsanov's theorem introduced in \cite{Schuppen1974}.

\proof[The proof of Lemma \ref{lemma 4.1 LLNunderTransformedMeasure}]

The existence and uniqueness of Equation \eqref{equ 4.1 measureValuedODE} is given in Appendix \ref{Appendix A.3}. We further prove in Appendix \ref{Appendix A.4} that $\{\theta^N\}_{N\geq 1}$ are $\hat{P}^N_{f, G}$-tight. Since $C^1([0, 1])$ is dense in $C([0, 1])$, we only need to check that if $\varpi$ is a $\hat{P}^N_{f, G}$-weak limit of a subsequence of $\{\theta^N\}_{N\geq 1}$, then $\varpi$ satisfies Equation \eqref{equ 4.1 measureValuedODE} for all $h\in C^1([0, 1])$.

According to It\^{o}'s formula and the definition of $\Gamma^N_{T_0}(G)$,
\begin{align}\label{equ 4.4}
d\Gamma^N_{T_0}(G)&=\frac{1}{H_G^N(0, X_0^N)}\exp\left\{\int_0^t\frac{(\partial_u+\mathcal{L}_N)H_G^N(u, X_u^N)}{H_G^N(u, X_u^N)}du\right\}d\mathcal{M}_t(H_G^N)\notag\\
&=\Gamma^N_{T_0}(G)d\widetilde{\mathcal{M}}_t(H_G^N),
\end{align}
where
\[
\widetilde{\mathcal{M}}_t(H_G^N)=\int_0^t\frac{1}{H_G^N(u, X_u^N)}d\mathcal{M}_u(H_G^N).
\]
For any $h\in C^1([0, 1])$, let
\[
\widehat{\mathcal{M}}_t(\theta^N(h))=\mathcal{M}_t(\theta^N(h))-\langle\mathcal{M}(\theta^N(h)), \widetilde{\mathcal{M}}(H_G^N)\rangle_t,
\]
then according to Equation \eqref{equ 4.4} and Theorem 3.2 of \cite{Schuppen1974}, which is a generalized version of Girsanov's theorem,
$\left\{\widehat{\mathcal{M}}_t(\theta^N(h))\right\}_{0\leq t\leq T_0}$ is a local martingale under $\hat{P}^{N}_{f, G}$ for all $h\in C^1([0, 1])$ and
\[
\left[\widehat{\mathcal{M}}(\theta^N(h)),\widehat{\mathcal{M}}(\theta^N(h))\right]_t=\left[\mathcal{M}(\theta^N(h)), \mathcal{M}(\theta^N(h))\right]_t
\]
under both $P_f^N$ and $\hat{P}^N_{f, G}$. Since $\{X^N_t\}_{t\geq 0}$ is a pure jump process and $EX^N_t(i)$ is differentiable with $t$ for $1\leq i\leq N$,
\[
\left[\mathcal{M}(\theta^N(h)), \mathcal{M}(\theta^N(h))\right]_t=\sum_{0\leq u\leq t}\left(\theta_u^N(h)-\theta_{u-}^N(h)\right)^2.
\]
Hence, by Lemma \ref{lemma 4.2} and Doob's inequality, $\sup_{0\leq t\leq T_0}|\widehat{\mathcal{M}}_t(\theta^N(h))|=o_p(1)$ under $\hat{P}^N_{f, G}$. By Equation \eqref{equ 4.3} and the definition of $\widetilde{\mathcal{M}}_t(H_G^N)$,
\begin{align*}
&d\langle\mathcal{M}(\theta^N(h)), \widetilde{\mathcal{M}}(H_G^N)\rangle_t\\
&=\sum_{i=1}^N\sum_{j=1}^N\frac{\lambda\left(\frac{i}{N}, \frac{j}{N}\right)X_t^N(i)}{N}\left(e^{\frac{a_N}{N}\left(G_t\left(\frac{j}{N}\right)-G_t\left(\frac{i}{N}\right)\right)}-1\right)\frac{h\left(\frac{j}{N}\right)
-h\left(\frac{i}{N}\right)}{a_N}dt.
\end{align*}
As a result, under $\hat{P}^N_{f, G}$,
\begin{align}\label{equ 4.5}
&\theta^N_t(h)=\theta^N_0(h)+o_p(1)+\int_0^t\theta_s^N\left((P_1^N-P_2^N)h\right)ds\\
&+\int_0^t\Bigg(\sum_{i=1}^N\sum_{j=1}^N\frac{\lambda\left(\frac{i}{N}, \frac{j}{N}\right)X_s^N(i)}{N}\Big(e^{\frac{a_N}{N}\left(G_s\left(\frac{j}{N}\right)-G_s\left(\frac{i}{N}\right)\right)}
-1\Big)\frac{h\left(\frac{j}{N}\right)
-h\left(\frac{i}{N}\right)}{a_N}\Bigg)ds.\notag
\end{align}
According to Taylor's expansion formula up to the second order,
\[
e^{\frac{a_N}{N}\left(G_s\left(\frac{j}{N}\right)-G_s\left(\frac{i}{N}\right)\right)}=\frac{a_N}{N}\left(G_s\left(\frac{j}{N}\right)-G_s\left(\frac{i}{N}\right)\right)
+O(\frac{a_N^2}{N^2}).
\]
Then, according to a similar analysis with that in the proof of Equation \eqref{equ 1.4 UpperBoundofMDP},
\begin{align*}
&\int_0^t\Bigg(\sum_{i=1}^N\sum_{j=1}^N\frac{\lambda\left(\frac{i}{N}, \frac{j}{N}\right)X_s^N(i)}{N}\Big(e^{\frac{a_N}{N}\left(G_s\left(\frac{j}{N}\right)-G_s\left(\frac{i}{N}\right)\right)}
-1\Big)\frac{h\left(\frac{j}{N}\right)
-h\left(\frac{i}{N}\right)}{a_N}\Bigg)ds\\
&+\int_0^t\theta_s^N((P_1^N-P_2^N)h)ds\\
&=\int_0^t\mu_s^N(\mathcal{R}_s(G,h))ds+\int_0^t\theta_s^N((P_1-P_2)h)ds+\epsilon_{7,t}^N,
\end{align*}
where $\mathcal{R}_s(G,h)(x)=\int_0^1\lambda(x, y)(G_s(y)-G_s(x))(h(y)-h(x))dy$ for all $0\leq x\leq 1$ and $\sup_{0\leq t\leq T_0}|\epsilon_{7,t}^N|=o_{\exp}\left(a_N\right)$ under $P$, the initial probability measure of our model. As we have shown in the proof of Lemma \ref{lemma 4.2}, conditioned on $\sum_{i=1}^NX_0^N(i)\leq NM$, $\frac{d\hat{P}^N_{f,G}}{dP^N_f}=\Gamma_{T_0}^N(G)\leq \exp\{a_NMC_6\}$ for some $C_6$ independent of $N$. Similarly, it is easy to check that
\begin{align*}
\frac{dP^N_f}{dP}&=\frac{e^{-\sum_{i=1}^N\left(\frac{a_N}{N}f\left(\frac{i}{N}\right)+\phi\left(\frac{i}{N}\right)\right)}
\prod_{i=1}^N\left(\frac{a_N}{N}f\left(\frac{i}{N}\right)+\phi\left(\frac{i}{N}\right)\right)^{X_0^N(i)}}
{e^{-\sum_{i=1}^N\phi\left(\frac{i}{N}\right)}\prod_{i=1}^N\phi\left(\frac{i}{N}\right)^{X_0^N(i)}}\\
&\leq \exp\{a_NMC_7\}
\end{align*}
for some $C_7=C_7(f)$ independent of $N$ conditioned on $\sum_{i=1}^NX_0^N(i)\leq NM$. Then it is easy to check that $\sup_{0\leq t\leq T_0}|\epsilon_{7,t}^N|=o_{\exp}\left(a_N\right)$ under $\hat{P}_{f, G}^{N}$ according to Equation \eqref{Appendix A.3} with $P^N_f$ replaced by $P$. By Lemma \ref{lemma 2.1 replacement},
\[
\int_0^t\mu_s^N(\mathcal{R}_s(G,h))ds=\int_0^t\mu_s(\mathcal{R}_s(G,h))ds+\epsilon_{8,t}^N,
\]
where $\sup_{0\leq t\leq T_0}|\epsilon_{8,t}^N|=o_{\exp}(a_N)$ under $P$. Then, according to a similar analysis with that of $\epsilon_{7,t}^N$, $\sup_{0\leq t\leq T_0}|\epsilon_{8,t}^N|=o_{\exp}\left(a_N\right)$ under $\hat{P}_{f, G}^N$. In conclusion, since $\mu_s(\mathcal{R}_s(G,h))=\langle G_s|h\rangle_s$,
\[
\theta_t^N(h)=\theta_0^N(h)+o_p(1)+\int_0^t\theta_s^N\left((P_1-P_2)h\right)ds+\int_0^t\langle G_s|h\rangle_sds
\]
under $\hat{P}^N_{f, G}$, where $o_p(1)$ can be chosen uniformly for $0\leq t\leq T_0$. Since we have proved that $\{\theta^N\}_{N\geq 1}$ is $\hat{P}_{f, G}^N$-tight in Appendix \ref{Appendix A.4}, we only need to show that
\begin{equation}\label{equ 4.6}
\theta_0^N(h)=\int_0^1f(x)h(x)dx+o_p(1)
\end{equation}
under $\hat{P}^N_{f, G}$ to finish the proof. As we have introduced in Section \ref{Section 3}, distributions of $\theta^N_0$ under $\hat{P}^N_{f, G}$ and $P^N_f$ are equal. As a result, Equation \eqref{equ 4.6} follows directly from the definition of $P^N_f$ and Chebyshev's inequality and hence the proof is complete.

\qed

At last we prove Equation \eqref{equ 1.5 LowerBoundofMDP}. The proof is a rigorous statement of the intuitive analysis given in Remark \ref{Remark 4.1}.

\proof[Proof of Equation \eqref{equ 1.5 LowerBoundofMDP}]

Equation \eqref{equ 1.5 LowerBoundofMDP} is trivial when $\inf_{\pi\in O}\{I_{ini}(\pi_0)+I_{dyn}(\pi)\}=+\infty$. Hence we only deal with the case where $\inf_{\pi\in O}\{I_{ini}(\pi_0)+I_{dyn}(\pi)\}<+\infty$. For any $\epsilon>0$, there exists $\pi^\epsilon\in O$ that
\[
I_{ini}(\pi^\epsilon_0)+I_{dyn}(\pi^\epsilon)\leq \inf_{\pi\in O}\{I_{ini}(\pi_0)+I_{dyn}(\pi)\}+\epsilon.
\]
Then, by Theorem \ref{theorem 1.3 alternative representation formula}, there exist $f^\epsilon\in L^2([0, 1])$ and $F^\epsilon\in \mathcal{H}$ that
$\pi_0^\epsilon(dx)=f^\epsilon(x)dx$,
\begin{equation}\label{equ 4.7}
\pi_{T_0}^\epsilon(G_{T_0})-\pi_0^\epsilon(G_0)-\int_0^{T_0}\pi_s^\epsilon\left((\partial_s+P_1-P_2)G_s\right)ds=\ll G, F^\epsilon\gg
\end{equation}
for any $G\in C^{1,0}([0, T_0]\times [0, 1])$ and
\[
I_{ini}(\pi_0^\epsilon)=\frac{1}{2}\int_0^1\frac{\left(f^\epsilon(x)\right)^2}{\phi(x)}dx, \text{~}I_{dyn}(\pi^\epsilon)=\frac{\ll F^\epsilon, F^\epsilon \gg}{2}.
\]
By Equation \eqref{equ 4.7}, let $G(s,x)=l(s)h(x)$ for some $l\in C^1([0, T_0])$ and $h\in C([0, 1])$, then
\[
l_{T_0}\pi^\epsilon_{T_0}(h)-l_0\pi^\epsilon_{0}(h)-\int_0^{T_0}l^\prime(s)\pi^\epsilon_s(h)ds=\int_0^{T_0}l(s)\pi^\epsilon_s((P_1-P_2)h)ds
+\int_0^{T_0}l(s)\langle F^\epsilon_s|h\rangle_sds.
\]
Since $l$ can be chosen arbitrarily, $\pi^\epsilon_t(h)$ is absolutely continuous with respect to $t$ and
\[
\frac{d}{dt}\pi^\epsilon_t(h)=\pi_t^\epsilon((P_1-P_2)h)+\langle F^\epsilon_t|h\rangle_t.
\]
Therefore, as we have shown in Appendix \ref{Appendix A.3},
\begin{equation}\label{equ 4.8}
\pi_t^\epsilon=e^{t(P_1-P_2)^{*}}\pi_0^\epsilon+\int_0^te^{(t-u)(P_1-P_2)^*}\Xi_u^{F^\epsilon}du.
\end{equation}
Since $C([0, 1])$ is dense in $L^2([0, 1])$ and $C^{1,0}([0, T_0]\times [0, 1])$ is dense in $\mathcal{H}$, there exist a sequence $\{f^n\}_{n\geq 1}$ in $C([0, 1])$ and a sequence $\{F^n\}_{n\geq 1}$ in $C^{1,0}([0, T_0]\times [0, 1])$ that $f^n\rightarrow f^\epsilon$ in $L^2$ and $F^n\rightarrow F^\epsilon$ in $\mathcal{H}$ and then,
\[
\lim_{n\rightarrow+\infty}\frac{1}{2}\int_0^1\frac{\left(f^n(x)\right)^2}{\phi(x)}dx=\frac{1}{2}\int_0^1 \frac{\left(f^\epsilon(x)\right)^2}{\phi(x)}dx, \text{~}\lim_{n\rightarrow+\infty}\frac{\ll F^n, F^n\gg}{2}=\frac{\ll F^\epsilon, F^\epsilon\gg}{2}.
\]
For each $n\geq 1$, let $\pi^n$ be the unique element in $\mathcal{D}([0, T_0], \mathcal{S})$ that
\begin{equation*}
\begin{cases}
\frac{d}{dt}\pi_t^n(h)=\pi_t^n\left((P_1-P_2)h\right)+\langle F_t^n | h\rangle_t \text{~for all~}0\leq t\leq T_0\text{~and~}h\in C([0, 1]),\\
\pi_0^n(dx)=f^n(x)dx.
\end{cases}
\end{equation*}
i.e., $\pi_t^n=e^{t(P_1-P_2)^{*}}\pi_0^n+\int_0^te^{(t-u)(P_1-P_2)^*}\Xi_u^{F^n}du$ as we have shown in Appendix \ref{Appendix A.3}. Then, by Equation \eqref{equ 4.8}, $\pi^n\rightarrow \pi^\epsilon$ in $\mathcal{D}([0, T_0], \mathcal{S})$. Furthermore, by Theorem \ref{theorem 1.3 alternative representation formula}, $I_{ini}(\pi^n_0)=\frac{1}{2}\int_0^1\frac{\left(f^n(x)\right)^2}{\phi(x)}dx$ and $I_{dyn}(\pi^n)=\frac{\ll F^n, F^n\gg}{2}$. Then, since $O$ is open, there exists $m\geq 1$ that $\pi^m\in O$ and
\[
I_{ini}(\pi^m_0)+I_{dyn}(\pi^m)\leq \inf_{\pi\in O}\{I_{ini}(\pi_0)+I_{dyn}(\pi)\}+2\epsilon.
\]
By Equation \eqref{equ 3.0},
\[
\Gamma_{T_0}^N(F^m)=\exp\left\{\frac{a_N^2}{N}\left(l(\theta^N, F^m)+\epsilon^N\right)\right\},
\]
where $\epsilon^N=o_{\exp}(a_N)$ under $P$. According to a similar analysis with that of $\epsilon_{7,t}^N$, it is easy to check that $\epsilon^N=o_{\exp}(a_N)$ under $\hat{P}^N_{f^m, F^m}$. Let
\[
D^\epsilon=\left\{\pi:~|l(\pi, F^m)-l(\pi^m, F^m)|<\epsilon\right\}\cap O,
\]
then $\hat{P}^N_{f^m, F^m}(\theta^N\in D^\epsilon)=1+o(1)$ as $N\rightarrow+\infty$ by Lemma \ref{lemma 4.1 LLNunderTransformedMeasure} and the fact that $\pi^m\in D^\epsilon$. According to the definition of $F^m$ and $\pi^m$, it is easy to check that $l(\pi^m, F^m)=\frac{\ll F^m, F^m\gg}{2}=I_{dyn}(\pi^m)$. By Chebyshev's inequality and the definition of $P^N_f$, it is easy to check that
\[
\frac{dP}{dP_{f^m}^N}=\exp\left\{-\frac{a_N^2}{N}\left(\frac{1}{2}\int_0^1\frac{\left(f^m(x)\right)^2}{\phi(x)}dx+o_p(1)\right)\right\}
\]
under $\hat{P}^{N}_{f^m, F^m}$. As a result, let
\[
\hat{D}^{\epsilon,N}=\{\theta^N\in D^\epsilon\}\cap\{|\epsilon^N|<\epsilon\}
\cap\left\{\frac{dP}{dP_{f^m}^N}\geq\exp\left\{-\frac{a_N^2}{N}\left(\frac{1}{2}\int_0^1\frac{\left(f^m(x)\right)^2}{\phi(x)}dx+\epsilon\right)\right\}\right\},
\]
then $\hat{P}^N_{f^m, F^m}(\hat{D}^{\epsilon,N})=1+o(1)$ as $N\rightarrow+\infty$ and
\[
\frac{dP}{d\hat{P}^N_{f^m, F^m}}=\left(\Gamma_{T_0}^N(F^m)\right)^{-1}\frac{dP}{dP_{f^m}^N}
\geq \exp\left\{-\frac{a_N^2}{N}\left(I_{ini}(\pi_0^m)+I_{dyn}(\pi^m)+3\epsilon\right)\right\}
\]
on $\hat{D}^{\epsilon, N}$. Since $\hat{D}^{\epsilon, N}\subseteq \{\theta^N\in O\}$,
\begin{align*}
P(\theta^N\in O)&\geq P(\hat{D}^{\epsilon, N})=E_{\hat{P}^N_{f^m, F^m}}\left(\frac{dP}{d\hat{P}^N_{f^m, F^m}}1_{\{\hat{D}^{\epsilon, N}\}}\right)\\
&\geq \exp\left\{-\frac{a_N^2}{N}\left(I_{ini}(\pi_0^m)+I_{dyn}(\pi^m)+3\epsilon\right)\right\}(1+o(1))
\end{align*}
and hence
\begin{align*}
\liminf_{N\rightarrow+\infty}\frac{N}{a_N^2}\log P(\theta^N\in O)&\geq -\left(I_{ini}(\pi_0^m)+I_{dyn}(\pi^m)\right)-3\epsilon\\
&\geq -\inf_{\pi\in O}\left(I_{ini}(\pi_0)+I_{dyn}(\pi)\right)-5\epsilon.
\end{align*}
Since $\epsilon$ is arbitrary, the proof is complete.

\qed

\appendix{}
\section{Appendix}
\subsection{Proof of Lemma \ref{lemma 3.1 independent and Poisson}}\label{subsection A.1}
 In this subsection we prove Lemma \ref{lemma 3.1 independent and Poisson}.

\proof[Proof of Lemma \ref{lemma 3.1 independent and Poisson}]

For $1\leq i, j\leq N$, we use $p_t^N(i,j)$ to denote probability that a given gas molecule is in the $j$th box at moment $t$ conditioned on it is in the $i$th box at moment $0$. Then, according to our assumption of $X_0^N$,
\begin{equation}\label{equ A.1}
EX_t^N(i)=\sum_{l=1}^NEX_0^N(l)p_t^N(l,i)=\sum_{l=1}^N\phi(\frac{l}{N})p_t^N(l,i).
\end{equation}
For $i\neq j$ and $1\leq k\leq X_0^N(i)$, we use $A_k^{N,t}(i, j)$ to denote the indicator function of the event that the $k$th gas molecule in the $i$th box at moment $0$ is in the $j$th box at moment $t$, then for given $r_1, r_2, \ldots, r_N\in \mathbb{R}$,
\[
\exp\left\{\sum_{j=1}^Nr_jX_t^N(j)\right\}=\exp\left\{\sum_{l=1}^N\sum_{k=1}^{X_0^N(l)}\sum_{j=1}^Nr_jA_k^{N,t}(l,j)\right\}.
\]
Therefore, according to our assumption of $X_0^N$ and Equation \eqref{equ A.1},
\begin{align*}
E\left(\exp\left\{\sum_{j=1}^Nr_jX_t^N(j)\right\}\Bigg|X_0^N\right)&=\prod_{l=1}^N\prod_{k=1}^{X_0^N(l)}\left(\sum_{j=1}^Ne^{r_j}p_t^N(l,j)\right) \\
&=\prod_{l=1}^N\left(\sum_{j=1}^Ne^{r_j}p_t^N(l, j)\right)^{X_0^N(l)}
\end{align*}
and
\begin{align}\label{equ A.2}
&E\left(\exp\left\{\sum_{j=1}^Nr_jX_t^N(j)\right\}\right)=\prod_{l=1}^NE\left(\left(\sum_{j=1}^Ne^{r_j}p_t^N(l, j)\right)^{X_0^N(l)}\right)\\
&=\prod_{l=1}^N \exp\left\{\left(\sum_{j=1}^Ne^{r_j}p_t^N(l, j)-1\right)\phi\left(\frac{l}{N}\right)\right\}
=\exp\left\{\sum_{j=1}^N\left(e^{r_j}-1\right)\sum_{l=1}^N\phi\left(\frac{l}{N}\right)p_t^N(l,j)\right\}\notag\\
&=\exp\left\{\sum_{j=1}^N\left(e^{r_j}-1\right)EX_t^N(j)\right\}. \notag
\end{align}
Since $r_1, r_2,\ldots, r_N$ are arbitrary, Lemma \ref{lemma 3.1 independent and Poisson} follows from Equation \eqref{equ A.2} directly.

\qed

\subsection{Proof of Lemma \ref{lemma 4.2}}\label{subsection A.2}

In this subsection, we prove Lemma \ref{lemma 4.2}.

\proof[Proof of Lemma \ref{lemma 4.2}]

We first show that Equation \eqref{equ 4.2} holds under $P^N_f$. For any $M>0$, according to Markov's inequality,
\[
\limsup_{N\rightarrow+\infty}\frac{1}{N}\log P^N_f\left(\frac{1}{N}\sum_{i=1}^NX_0^N(i)\geq M\right)\leq (e-1)\int_0^1\phi(x)dx-M
\]
and hence
\begin{equation}\label{equ A.3}
\limsup_{M\rightarrow+\infty}\limsup_{N\rightarrow+\infty}\frac{1}{N}\log P^N_f\left(\frac{1}{N}\sum_{i=1}^NX_0^N(i)\geq M\right)=-\infty.
\end{equation}
Note that Equation \eqref{Appendix A.3} still holds when $P^N_f$ is replaced by the original probability measure $P$ of our process according to the same analysis as that under $P^N_f$. Conditioned on $\frac{1}{N}\sum_{i=1}^NX_0^N(i)\leq M$, $\theta^N_t$ jumps at rate at most $\|\lambda\|NM$, where $\|\lambda\|=\sup_{0\leq x,y\leq 1}|\lambda(x,y)|$ and
\[
\left(\theta_t^N(h)-\theta_{t-}^N(h)\right)^2\leq \frac{4}{a^2_N}\|h\|^2
\]
when $t$ is a jump moment, where $\|h\|=\sup_{0\leq x\leq 1}|h(x)|$. As a result, conditioned on $\frac{1}{N}\sum_{i=1}^NX_0^N(i)\leq M$, $\sum_{0\leq t\leq T_0}\left(\theta_t^N(h)-\theta_{t-}^N(h)\right)^2$ is stochastically dominated from above by $\frac{4\|h\|^2}{a_N^2}Y(NM\|\lambda\|T_0)$, where $\{Y(t)\}_{t\geq 0}$ is a Poisson process with rate $1$. Hence, by Markov's inequality,
\begin{align*}
&P_f^N\left(\sum_{0\leq t\leq T_0}\left(\theta_t^N(h)-\theta_{t-}^N(h)\right)^2\geq \epsilon, \sum_{i=1}^NX_0^N(i)\leq NM\right)\leq P\left(Y(NM\|\lambda\|T_0)\geq \frac{a_N^2\epsilon}{4\|h\|^2}\right) \\
&\leq e^{-\frac{a_N^2\epsilon}{4\|h\|^2}}e^{(e-1)NMT_0\|\lambda\|}
\end{align*}
and consequently
\[
\limsup_{N\rightarrow+\infty}\frac{1}{N}\log P_f^N\left(\sum_{0\leq t\leq T_0}\left(\theta_t^N(h)-\theta_{t-}^N(h)\right)^2\geq \epsilon, \sum_{i=1}^NX_0^N(i)\leq NM\right)=-\infty
\]
since $\frac{a_N^2}{N}\rightarrow+\infty$. As a result,
\begin{align*}
&\limsup_{N\rightarrow+\infty}\frac{1}{N}\log P_f^N\left(\sum_{0\leq t\leq T_0}\left(\theta_t^N(h)-\theta_{t-}^N(h)\right)^2\geq \epsilon\right)\\
&\leq \limsup_{N\rightarrow+\infty}\frac{1}{N}\log P_f^N\left(\frac{1}{N}\sum_{i=1}^NX_0^N(i)\geq M\right).
\end{align*}
Since $M$ is arbitrary, Equation \eqref{equ 4.2} holds under $P_f^N$ according to Equation \eqref{equ A.3}.

\quad

Now we prove that Equation \eqref{equ 4.2} holds under $\hat{P}^N_{f, G}$. Conditioned on $\frac{1}{N}\sum_{i=1}^NX_0^N(i)\leq M$, it is easy to check that there exists $C_6=C_6(G)<+\infty$ independent of $N$ that $\Gamma_{T_0}^N(G)\leq e^{a_NC_6M}$ for sufficiently large $N$. Then, for sufficiently large $N$,
\begin{align*}
&\hat{P}_{f, G}^N\left(\sum_{0\leq t\leq T_0}\left(\theta_t^N(h)-\theta_{t-}^N(h)\right)^2\geq \epsilon, \frac{1}{N}\sum_{i=1}^NX_0^N(i)\leq M\right) \\
&=E_{P_f^N}\left(\Gamma_{T_0}^N(G)1_{\left\{\sum_{0\leq t\leq T_0}\left(\theta_t^N(h)-\theta_{t-}^N(h)\right)^2\geq \epsilon, \frac{1}{N}\sum_{i=1}^NX_0^N(i)\leq M\right\}}\right)\\
&\leq e^{a_NC_6M}P_f^N\left(\sum_{0\leq t\leq T_0}\left(\theta_t^N(h)-\theta_{t-}^N(h)\right)^2\geq \epsilon\right).
\end{align*}
Since we have shown that Equation \eqref{equ 4.2} holds under $P_f^N$ and $\lim_{N\rightarrow+\infty}\frac{a_N}{N}=0$,
\[
\lim_{N\rightarrow+\infty} \frac{1}{N}\log \hat{P}_{f, G}^N\left(\sum_{0\leq t\leq T_0}\left(\theta_t^N(h)-\theta_{t-}^N(h)\right)^2\geq \epsilon, \frac{1}{N}\sum_{i=1}^NX_0^N(i)\leq M\right)=-\infty
\]
and hence
\begin{align*}
&\limsup_{N\rightarrow+\infty}\frac{1}{N}\log \hat{P}_{f,G}^N\left(\sum_{0\leq t\leq T_0}\left(\theta_t^N(h)-\theta_{t-}^N(h)\right)^2\geq \epsilon\right)\\
&\leq \limsup_{N\rightarrow+\infty}\frac{1}{N}\log \hat{P}_{f,G}^N\left(\frac{1}{N}\sum_{i=1}^NX_0^N(i)\geq M\right).
\end{align*}
Since $\Gamma_0^N(G)=1$ and $\{\Gamma_t^N(G)\}_{0\leq t\leq T_0}$ is a martingale, distributions of $X_0^N$ under $P_f^N$ and $\hat{P}_{f, G}^N$ are the same. Then, since $M$ is arbitrary, Equation \eqref{equ 4.2} holds under $\hat{P}_{f, G}^N$ according to Equation \eqref{equ A.3}.

\qed

\subsection{Existence and uniqueness of the solution to Equation \eqref{equ 4.1 measureValuedODE}}\label{Appendix A.3}

In this subsection we give the proof of existence and uniqueness of the solution to Equation \eqref{equ 4.1 measureValuedODE}.

\proof[Proof of existence and uniqueness of the solution to Equation \eqref{equ 4.1 measureValuedODE}] For any $\mu \in \mathcal{S}$, we use $\|\mu\|$ to denote the norm of $\mu$, i.e.,
\[
\|\mu\|=\sup\left\{|\mu(f)|:~f\in C([0, 1])\text{~and~}\sup_{0\leq x\leq 1}|f(x)|\leq 1\right\}.
\]
We further define $(P_1-P_2)^{*}$ as the linear operator from $\mathcal{S}$ to $\mathcal{S}$ that
\[
\left((P_1-P_2)^{*}\mu\right)(f)=\mu\left((P_1-P_2)f\right)
\]
for any $\mu\in \mathcal{S}$ and $f\in C([0, 1])$. Then it is easy to check that $\|(P_1-P_2)^{*}\mu\|\leq 2\|\lambda\|\|\mu\|$ for any $\mu\in \mathcal{S}$. As a result, it is reasonable to define
\[
e^{c(P_1-P_2)^{*}}=\sum_{n=0}^{+\infty}\frac{c^n((P_1-P_2)^{*})^n}{n!}
\]
for any $c\in \mathbb{R}$ and the domain of $e^{c(P_1-P_2)^{*}}$ is $\mathcal{S}$.
For $G\in C^{1,1}([0, T_0]\times[0, 1])$ and any $0\leq t\leq T_0$, let $\Xi_t^G$ be the element in $\mathcal{S}$ that $\Xi_t^G(f)=\langle G_t|f\rangle_t$
for any $f\in C([0, 1])$. Then Equation \eqref{equ 4.1 measureValuedODE} can be considered as a $\mathcal{S}$-valued linear ODE that
\[
\begin{cases}
&\frac{d}{dt}\vartheta^{f,G}_t=(P_1-P_2)^{*}\vartheta^{f,G}_t+\Xi_t^G\text{~for~}0\leq t\leq T_0,\\
&\vartheta^{f,G}_0(dx)=f(x)dx.
\end{cases}
\]
Therefore,
\[
\frac{d}{dt}\left(e^{-t(P_1-P_2)^{*}}\vartheta^{f,G}_t\right)=e^{-t(P_1-P_2)^{*}}\Xi_t^G
\]
and hence
\[
\vartheta_t^{f,G}=e^{t(P_1-P_2)^{*}}\vartheta_0^{f,G}+\int_0^te^{(t-u)(P_1-P_2)^{*}}\Xi_u^Gdu,
\]
where $\vartheta^{f,G}_0(dx)=f(x)dx$. Since we have directly solved Equation \eqref{equ 4.1 measureValuedODE}, the solution exists and is unique.

\qed

\subsection{$\hat{P}^N_{f, G}$-tightness of $\{\theta^N\}_{N\geq 1}$}\label{Appendix A.4}

In this subsection we prove that $\{\theta^N\}_{N\geq 1}$ is $\hat{P}^N_{f, G}$-tight.

\proof[Proof of $\hat{P}^N_{f, G}$-tightness of $\{\theta^N\}_{N\geq 1}$]

By Aldous' criteria, we only need to check that the following two claims hold.

\textbf{Claim} 1. For all $h\in C([0, 1])$,
\[
\lim_{M\rightarrow+\infty}\limsup_{N\rightarrow+\infty} \hat{P}^N_{f, G}\left(|\theta^N_t(h)|\geq M\right)=0
\]
for all $0\leq t\leq T_0$.

\textbf{Claim} 2. For any $\epsilon>0$ and $h\in C([0, 1])$,
\[
\lim_{\delta\rightarrow 0}\limsup_{N\rightarrow+\infty}\sup_{\tau\in \Upsilon, s\leq \delta}\hat{P}^N_{f, G}\left(|\theta^N_{\tau+s}(h)-\theta^N_\tau(h)|>\epsilon\right)=0,
\]
where $\Upsilon$ is the set of stopping times of $\{X_t^N\}_{t\geq 0}$ bounded by $T_0$.

We first check Claim 1. As we have shown in Sections \ref{Section 3} and \ref{section 4},
\[
\Gamma_{T_0}^N(G)=\exp\left\{\frac{a_N^2}{N}\left(l(\theta^N, G)+\epsilon^N\right)\right\},
\]
where $\epsilon^N=o_{\exp}(a_N)$ under both $P$ and $\hat{P}^N_{f, G}$. Hence, to check Claim 1, we only need to show that
\begin{equation}\label{equ A.4}
\lim_{M\rightarrow+\infty}\limsup_{N\rightarrow+\infty} \hat{P}^N_{f, G}\left(|\theta^N_t(h)|\geq M, |\epsilon^N|\leq 1\right)=0.
\end{equation}
By H\"{o}lder's inequality, Markov's inequality and the fact that
\[
\left(\frac{d\hat{P}^N_{f, G}}{dP^N_f}\right)^2=\left(\Gamma_{T_0}^N(G)\right)^2\leq\exp\left\{\frac{2a_N^2}{N}\left(l(\theta^N, G)+1\right)\right\}
\]
when $|\epsilon^N|\leq 1$, to prove Equation \eqref{Appendix A.4} we only need to show that
\begin{equation}\label{equ A.5}
\limsup_{N\rightarrow+\infty}\frac{N}{a_N^2}\log\sup_{0\leq t\leq T_0}E_{_{P^N_f}}\left(\exp\left\{\frac{a_N^2}{N}\theta_t^N(h)\right\}\right)<+\infty
\end{equation}
and
\begin{equation}\label{equ A.6}
\limsup_{N\rightarrow+\infty}\frac{N}{a_N^2}\log\sup_{0\leq t\leq T_0}E_{_{P^N_f}}\left(\exp\left\{\frac{Ca_N^2}{N}l(\theta^N, G)\right\}\right)<+\infty
\end{equation}
for any $C>0$.
By Lemma \ref{lemma 3.1 independent and Poisson}, under $P^N_f$, $\{X_t^N(i)\}_{1\leq i\leq N}$ are independent and $X_t^N(i)$ follows Poisson distribution with mean
\[
E_{_{P^N_f}} X_t^N(i)=\sum_{j=1}^N\left(\phi\left(\frac{j}{N}\right)+\frac{a_N}{N}f\left(\frac{j}{N}\right)\right)p_t^N(j,i)
\]
for all $1\leq i\leq N$. As a result,
\begin{align*}
&E_{_{P^N_f}}\left(\exp\left\{\frac{a_N^2}{N}\theta_t^N(h)\right\}\right)\\
&=e^{\sum_{i=1}^NE_{_{P^N_f}}X_t^N(i)\left(e^{\frac{a_N}{N}h(\frac{i}{N})}-\frac{a_N}{N}h(\frac{i}{N})-1\right)
+\sum_{i=1}^N\frac{a_N}{N}h(\frac{i}{N})\left(E_{_{P^N_f}}X_t^N(i)-EX_t^N(i)\right)}.
\end{align*}
Since $\sum_{i=1}^Np_t^N(j,i)=1$,
\[
\sum_{i=1}^N\frac{a_N}{N}h(\frac{i}{N})\left(E_{_{P^N_f}}X_t^N(i)-EX_t^N(i)\right)\leq \frac{a_N^2}{N}\|h\|\|f\|.
\]
According to a similar analysis with that given in Section 4 of \cite{Xue2020}, there exists $C_9$ independent of $N$ such that
\[
\sup_{1\leq i\leq N, 0\leq t\leq T_0}E_{_{P^N_f}}X_t^N(i)\leq C_9
\]
for sufficiently large $N$. Therefore, according to Taylor's expansion formula up to the second order,
\[
\limsup_{N\rightarrow+\infty}\frac{N}{a_N^2}\log\sup_{0\leq t\leq T_0}E_{_{P^N_f}}\left(\exp\left\{\frac{a_N^2}{N}\theta_t^N(h)\right\}\right)\leq
C_9\int_0^1h^2(x)dx+\|h\|\|f\|
\]
and hence Equation \eqref{equ A.5} holds. Now we check Equation \eqref{equ A.6}. By repeated utilizing H\"{o}lder's inequality and Jensen's inequality,
\begin{align}\label{equ A.7}
&E_{_{P_f^N}}\left(\exp\left\{\frac{Ca_N^2}{N}\theta_{T_0}^N(G_{T_0})-\frac{Ca_N^2}{N}\theta_0^N(G_0)
-\int_0^{T_0}\frac{Ca_N^2}{N}\theta_s^N\left((\partial_s+P_1-P_2)G_s\right)ds\right\}\right) \notag\\
&\leq \sqrt{E_{_{P_f^N}}e^{\frac{a_N^2}{N}\theta_{T_0}^N(2CG_{T_0})+\frac{a_N^2}{N}\theta_0^N(-2CG_0)}}\sqrt{E_{_{P_f^N}}
e^{\int_0^{T_0}\frac{a_N^2}{N}\theta_s^N\left(-2C(\partial_s+P_1-P_2)G_s\right)ds}} \\
&\leq \left(E_{_{P_f^N}}e^{\frac{a_N^2}{N}\theta_{T_0}^N(4CG_{T_0})}\right)^{\frac{1}{4}}\left(E_{_{P_f^N}}e^{\frac{a_N^2}{N}\theta_0^N(-4CG_0)}\right)^{\frac{1}{4}}
\notag\\
&\text{\quad}\times\sqrt{E_{_{P_f^N}}\left(\frac{1}{T_0}\int_0^{T_0}e^{\frac{a_N^2}{N}\theta_s^N\left(-2CT_0(\partial_s+P_1-P_2)G_s\right)}ds\right)} \notag\\
&= \left(E_{_{P_f^N}}e^{\frac{a_N^2}{N}\theta_{T_0}^N(4CG_{T_0})}\right)^{\frac{1}{4}}\left(E_{_{P_f^N}}e^{\frac{a_N^2}{N}\theta_0^N(-4CG_0)}\right)^{\frac{1}{4}}
\notag\\
&\text{\quad}\times\sqrt{\frac{1}{T_0}\int_0^{T_0}\left(E_{_{P_f^N}}e^{\frac{a_N^2}{N}\theta_s^N\left(-2CT_0(\partial_s+P_1-P_2)G_s\right)}\right)ds}. \notag
\end{align}
Equation \eqref{equ A.6} follows from Equations \eqref{equ A.5} and \eqref{equ A.7} and hence Claim 1 holds.

Now we check Claim 2. As we have shown in Section \ref{section 4}, under $\hat{P}^N_{f, G}$,
\[
\theta^N_t(h)=\int_0^1f(x)h(x)dx+o_p(1)+\int_0^t\theta_s^N((P_1-P_2)h)ds+\int_0^t\langle G_s|h\rangle_s ds,
\]
where $o_p(1)$ can be chosen uniformly for $0\leq t\leq T_0$. Hence, to prove Claim 2, we only need to check that
\begin{equation}\label{equ A.8}
\lim_{\delta\rightarrow 0}\limsup_{N\rightarrow+\infty}\sup_{\tau\in \Upsilon, s\leq \delta}\hat{P}^N_{f, G}\left(\left|\int_{\tau}^{\tau+s}\theta_u^N((P_1-P_2)h)du\right|>\epsilon\right)=0.
\end{equation}

As we have shown in Section \ref{section 4},
\[
\frac{dP}{dP^N_f}=\exp\left\{-\frac{a_N^2}{N}\left(\frac{1}{2}\int_0^1\frac{f^2(x)}{\phi(x)}dx+\epsilon_9^N\right)\right\},
\]
where $\epsilon_9^N=o_p(1)$ under $P^f_N$ and $\hat{P}^f_N$. Hence, to prove Equation \eqref{equ A.8}, we only need to check that
\begin{equation}\label{equ A.9}
\lim_{\delta\rightarrow 0}\limsup_{N\rightarrow+\infty}\sup_{\tau\in \Upsilon, s\leq \delta}\hat{P}^N_{f, G}\left(\left|\int_{\tau}^{\tau+s}\theta_u^N((P_1-P_2)h)du\right|>\epsilon, |\epsilon^N|\leq 1, |\epsilon_9^N|\leq 1\right)=0.
\end{equation}
By H\"{o}lder's inequality,
\begin{align*}
&\hat{P}^N_{f, G}\left(\left|\int_{\tau}^{\tau+s}\theta_u^N((P_1-P_2)h)du\right|>\epsilon, |\epsilon^N|\leq 1, |\epsilon_9^N|\leq 1\right)\\
&\leq\sqrt{E_{_{P^N_f}}e^{\frac{2a_N^2}{N}\left(l(\theta^N, G)+1\right)}}\sqrt{P_f^N\left(\left|\int_{\tau}^{\tau+s}\theta_u^N((P_1-P_2)h)du\right|>\epsilon, |\epsilon^N_9|\leq 1\right)}\\
&\leq \sqrt{E_{_{P^N_f}}e^{\frac{2a_N^2}{N}\left(l(\theta^N, G)+1\right)}}\sqrt{e^{\frac{a_N^2}{N}\left(1+\frac{1}{2}\int_0^1\frac{f^2(x)}{\phi(x)}dx\right)}}
\sqrt{P\left(\left|\int_{\tau}^{\tau+s}\theta_u^N((P_1-P_2)h)du\right|>\epsilon\right)}\\
&\leq \sqrt{E_{_{P^N_f}}e^{\frac{2a_N^2}{N}\left(l(\theta^N, G)+1\right)}}\sqrt{e^{\frac{a_N^2}{N}\left(1+\frac{1}{2}\int_0^1\frac{f^2(x)}{\phi(x)}dx\right)}}\\
&\text{\quad}\times \sqrt{P\left(\sup_{0\leq t_1<t_2\leq T_0, \atop |t_2-t_1|<\delta}\left|\int_{t_1}^{t_2}\theta_u^N((P_1-P_2)h)du\right|>\epsilon\right)}.
\end{align*}
As a result, Equation \eqref{equ A.9} follows from Lemma \ref{lemma 3.2 control of integration} and Equation \eqref{equ A.6} and hence Claim 2 holds. Since Claims 1 and 2 both hold, the proof is complete. 

\qed

\quad

\textbf{Acknowledgments.} The author is grateful to the financial support from the National Natural Science Foundation of China with grant number 11501542.

{}
\end{document}